\documentclass[11pt]{amsart}

\usepackage{amsmath}
\usepackage{amsfonts}
\usepackage{amssymb}

\def \End {\mathop {\rm End}}
\def \Hom {\mathop{\rm Hom}\nolimits}
\def \Homod {\Hom^{\rm odd}} 
\def \Asso {\mathop {\rm Asso}}
\def \Set {\mathop {\rm Set}}
\def \Com{\mathop{\rm Com}}
\def \Mod{\mathop{\rm Mod}}
\def \Lie{\mathop{\rm Lie}}
\def \Unit{\mathop{\rm Unit}}

\def \im {\mathop {\rm im}}

\newtheorem{theorem}{Theorem}[section]

\newtheorem{proposition}[theorem]{Proposition}
\newtheorem{corollary}[theorem]{Corollary}
\newtheorem{definition}[theorem]{Definition\rm}

\begin{document}
\title[Equivariant cohomology]{Equivariant cohomology \\ and tensor categories}
\author{Martin Andler}
\address{D\'{e}partement de Math\'ematiques (UMR CNRS 8100),
Universit\'{e} de Versailles Saint-Quentin, 78035
Versailles C\'{e}dex}
\email{andler@math.uvsq.fr}
\thanks{M. A. wishes to thank Rutgers University and the NSF for an invitation during which part of this research was conducted }
\author{Siddhartha Sahi}
\address{Mathematics Department, Rutgers
University, New Brunswick, NJ 08903}
\email{sahi@math.rutgers.edu}
\thanks{The research of S. S. was supported by an NSF grant; he wishes to thank Universit\'e de Versailles Saint-Quentin for an invitation during which part of this research was conducted}

\begin{abstract}
In this research announcement we propose the notion of a \emph{supercategory} as an
alternative approach to supermathematics. 
We show that this setting is rich to carry out 
many of the basic constructions of supermathematics. We also  prove generalizations of a number of results in
equivariant cohomology, including the Chern-Weil theorem for an arbitrary
rigid Lie algebra object. For a quadratic Lie algebra object we
obtain a proof of the Duflo isomorphism along the lines of Alekseev-Meinrenken, thereby
generalizing their  result to Lie superalgebras.
\end{abstract}
\date{February 7, 2008}
\maketitle

\section{Introduction}
\label{}
As pointed out in \cite{demo} many constructions of multilinear algebra
and differential geometry can be carried out in the context of a tensor
category, and the more general perspective often yields additional insights.
One advantage of this approach is that it treats vector spaces and super
vector spaces on an equal footing, thus avoiding
the signs which proliferate in the former setting. However one sometimes needs the
signs in order to deal with essentially ``super'' concepts, such as that of a connection.

In this research announcement we propose the notion of a \emph{supercategory} as an
alternative approach to super-mathematics. We show that this setting is rich
enough to carry out many of the basic constructions; and to further
demonstrate its power we establish a number of results in
equivariant cohomology, including the Chern-Weil theorem for an arbitrary
rigid Lie algebra object. In the case of a quadratic Lie algebra object we obtain
a proof of the Duflo isomorphism (see \cite{dufl}) along the lines of \cite{am2}. Our result implies in particular a generalization of the Alekseev-Meinrenken result to the setting of Lie superalgebras. The Duflo isomorphism for Lie algebras has been generalized in the Kashiwara-Vergne Conjecture (now Theorem~: see \cite{kash-ver}, \cite{ver}, \cite{to} and the references therein). It would be interesting to see whether this admits an extension to supercategories.

The proofs will appear elsewhere.

\subsection{Supercategories}

A tensor category is an abelian category equipped with a bilinear product
$\otimes,$ which is associative, commutative and has a unit $I$. We will
write $XY$ for $X\otimes Y$, and $X(\_)  $ for the functor
$Y\mapsto$ $XY$. A tensor category is said to be $\mathbb{Q}$-linear if there
is an inclusion $\mathbb{Q} \subset \End( I)  $ (this implies that
all the hom-sets are $\mathbb{Q}$-vector spaces). For basic facts about tensor
categories, including definitions of tensor, symmetric and exterior powers
$T^{\ast}(  X)  ,S^{\ast}(  X)  ,\Lambda^{\ast}(X)  $, we refer the reader to \cite{demo}.

\begin{definition}
A supercategory is a well-powered and complete $\mathbb{Q}$-linear tensor
category $\mathcal{S}$, equipped with the choice of an object $P$ satisfying
$\Lambda^{2}(  P)  \approx I,S^{2}(  P)  \approx 0$.
\end{definition}

As developed in \cite[p.45]{demo}, an important example of a supercategory is the category of superspaces over a field $k$ of characteristic $0.$ The
category admits a \emph{parity reversal functor} $\Pi$.
Objects in a general supercategory need not be bigraded. However we still have
a parity reversal functor, viz. $P(\cdot)  $, which satisfies
$P^{2}X\approx X$. One can also define the concept of an \emph{odd morphism}
from $X$ to $Y$, an element of
\begin{equation}\label{=homod}
\Homod (  X,Y)  :=\Hom(  PX,Y)  \approx \Hom(X,PY)
\end{equation}
where the two ordinary hom-spaces are identified via $P$. For $f\in
\Homod (X,Y)  $ and $g\in \Homod (Y,Z)$, their
composition $gf$ is an ordinary morphism in
$\Hom(  X,Y)  \approx \Hom(  PX,PY)$.

\subsection{Algebras and modules}

An ``ordinary'' algebra or module consists of a vector space over a field $k$,
with some algebraic structure which
can be formulated in terms of linear maps. Thus the definitions continue to make sense if we
replace the category of $k$-vector spaces by a tensor category. In particular
for a supercategory $\mathcal{S}$ we can define categories $\Asso(
\mathcal{S}), \Com(  \mathcal{S}), \Lie(\mathcal{S})$ of associative, commutative and Lie algebras,
respectively. 

Objects in these algebra categories have structure morphisms: 1)
product/Lie bracket: $A\otimes A\rightarrow A$ and 2) (for $\Asso(  \mathcal{S}), \Com(  \mathcal{S})$) unit: $I\rightarrow A$. For
each algebra $A$, one can also define a category $\Mod(A)$ of
$A$-modules; these consist of an underlying object $M$ in $\mathcal{S}$ and a
morphism $A\otimes M\overset{a}{\rightarrow}M$ called ``action'', required to satisfy the usual conditions, formulated in categorical terms.

\begin{proposition}
If $L\in \Lie(\mathcal{S})  $ then $\Mod(L)$ is a supercategory with a forgetful functor $F_{L}$
to $\mathcal{S}$. There is a trivial extension functor $E_{L}:\mathcal{S}\rightarrow \Mod(L)  $ where for every $X$ in
$\mathcal{S}$ the $L$-action on $E_{L}(  X)  $ is $0$. Furthermore, $E_{L}$ admits a right adjoint $M\mapsto
M^{L}: \Mod(L)  \rightarrow\mathcal{S}$, called the
``invariants functor''.
\end{proposition}

This proposition underlies an important philosophical point in our approach :  $L$-equivariant analogs of various algebraic constructions can be regarded as ``ordinary'' constructions in the supercategory $\Mod(L)$, and this point of view leads automatically to the correct definitions.

\section{Cohomology}

\subsection{The category $\mathcal{Q}$ and cohomology}

We fix a Lie algebra $Q$ whose underlying object is isomorphic to $P$, and
which is abelian, i.e., the bracket $Q\times Q\rightarrow Q$ is the $0$ morphism. One sees easily that
the supercategory $\mathcal{Q}= \Mod(Q)  $ is naturally isomorphic
to the category of \ ``supercomplexes''. It consists of pairs $(
X,d)  $, where $X$ is an object in $\mathcal{S}$ and $d$ (the
``differential'') is an odd endomorphism $d\in$ $\Homod (X,X)
$ satisfying $d^{2}=0$. We may therefore define a functor $\mathcal{Q}\rightarrow\mathcal{S}$, the $\mathcal{Q}$-\emph{cohomology functor}, as follows:
\begin{equation}\label{=hom}
H^{\mathcal{Q}}(X)  =(\ker d)/(\im d).
\end{equation}

Since the second and higher symmetric powers of $P$ vanish, the symmetric
algebra $D=S(  P)  $ has underlying object $I\oplus P$. In fact $D$
has a natural commutative $\mathcal{Q}$-algebra structure with differential
$d$ given as follows: $I\overset{d}{\rightarrow}0$, while $P\overset
{d}{\rightarrow}I$ is the odd endomorphism corresponding to the identity.$\ $%
Tensoring with $D$ defines the \emph{doubling functor} from $\mathcal{S}$ to
$\mathcal{Q}$, which takes $\mathcal{S}$-algebras to $\mathcal{Q}$-algebras
and $A$-modules to $DA$-modules.

\subsection{$\mathcal{L}$-cohomology and homotopy}

For the rest of the note we fix a Lie algebra $L_{0}\in \Lie(
\mathcal{S})  $. We write $L$ for $DL_{0}\in \Lie(  \mathcal{Q}
)  $ and $\mathcal{L}$ for the supercategory $\Mod(  L)  $. We
define the $\mathcal{L}$-cohomology functor $H^{\mathcal{L}}(  X)
:=H^{\mathcal{Q}}(  X^{L})  $ from $\mathcal{L}\rightarrow {S}$, and
for any $Z\in\mathcal{L}$ we define the twisted cohomology functor
$H_{Z}^{\mathcal{L}}(  X)  :=H^{\mathcal{L}}(  ZX)  .$
If $L_{0}$ is understood we will simply write $H$ and $H_{Z}$ for the two functors.

A \emph{homotopy} between two $\mathcal{L}$-morphisms $f,g:X\rightarrow Y$ is an odd
$\mathcal{S}$-morphism $h:X\rightarrow Y$ which commutes with the $L$-action
and satisfies $dh+hd=f-g$. An $\mathcal{L}$-morphism $f^{\prime
}:Y\rightarrow X$ is said to be the \emph{homotopy inverse} of $f:X\rightarrow Y$ if
$ff^{\prime}$ and $f^{\prime}f$ are homotopic to the identity morphisms for
$Y$ and $X$ respectively.

\begin{proposition}
Suppose $f:A\rightarrow A^{\prime}$ is an algebra morphism in $\Asso(
\mathcal{L})$, $B$ is an algebra in $\Asso(\mathcal{L})$, and $C$\ is a $B$-module; then
\begin{enumerate}
\item $H_{A}(  B)  $ is an algebra in $\Asso(\mathcal{S})$ and $H_{A}(C)$ is an $H_{A}(B)$
module. 
\item $f$ induces an algebra morphism $\tau_{B}(  f)  :H_{A}(
B)  \rightarrow H_{A^{\prime}}(  B)  $ and a compatible
module morphism $\tau_{C}(f)  : H_{A}(C)  \rightarrow H_{A^{\prime}}(C)  .$

\item  If $g$ is homotopic to $f$ then $\tau_{B}(  f)  =\tau_{B}(  g)  $ and $\tau_{C}(  f)  =\tau_{C}(g).$
\end{enumerate}
\end{proposition}

\subsection{Quasi-isomorphism}

An algebra morphism in $\Asso(  \mathcal{L})  $ is said to be a
\emph{quasi-isomorphism} if it admits a homotopy inverse as an $\mathcal{L}$-morphism. As it turns out, many algebra morphisms defined by universal properties turn out to be
quasi-isomorphisms, emphasizing the importance of : 
\begin{proposition}
\label{quasi} Suppose $f:A\rightarrow A^{\prime}$ is a quasi-isomorphism.

\begin{enumerate}
\item  If $B$ is an algebra in $\mathcal{L}$ then $\tau_{B}(  f)
:H_{A}(  B)  \rightarrow H_{A^{\prime}}(  B)  $ is an
algebra isomorphism.
\item  If $C$ is a $B$-module, then $\tau_{C}(  f)  :H_{A}(
C)  \rightarrow H_{A^{\prime}}(  C)  $ is a module isomorphism.
\end{enumerate}
\end{proposition}

\section{Chern-Weil theory}
\subsection{Connections}
Suppose $G$ is a compact Lie group with Lie algebra $\mathfrak{g}$ and $E$ is a
principal $G$-bundle. A (principal) connection on $E$ is a $G$-map $\theta
$ from $\mathfrak{g}^{\ast}$to $\Omega(  E)$, the de Rham complex of
differential forms on $E$, satisfying certain properties. This definition can
be generalized to the present situation, but in order to carry this out we
need the concepts of ``dual object'' and ``unital object''.

In a tensor categoy, an object $Y$ is said to be \emph{dual} to $X$ if there
exist compatible morphisms $I\rightarrow X\otimes Y$ (inclusion) and $Y\otimes
X\rightarrow I$ (evaluation). An object $X$ is \emph{rigid} if it has a dual, which is then unique. If so, $X$ is reflexive. For the category of vector spaces, rigidity corresponds to finite dimensionality.

A \emph{unital object} is an object $X$ in $\mathcal{S}$ together with a
distinguished morphism $u:I\rightarrow X$. There is an obvious notion of
unital morphisms between unital objects and hence unital objects in
$\mathcal{S}$ form a category $\Unit(  \mathcal{S})$. Moreover,
there is a natural sequence of forgetful functors $\Com(\mathcal{S})  \rightarrow \Asso(  \mathcal{S})  \rightarrow \Unit(\mathcal{S})  \rightarrow\mathcal{S}$.

\begin{definition}
Assume that $L_{0}$ is \emph{rigid}. A \emph{connection} on an object $X$ in $\Unit(\mathcal{L})
$ is a $\theta\in \Homod(L_{0}^{\ast},X)$ for
which the diagram commutes (the vertical maps are the evaluation and the $L$ action on $X$, respectively):
\begin{equation}
\begin{array}
[c]{ccc}
PL_{0}\otimes PL_{0}^{\ast} & \overset{1\otimes\theta}{\longrightarrow} &
PL_{0}\otimes X\\
\downarrow & \circlearrowleft & \downarrow\\
I & \overset{u}{\rightarrow} & X.
\end{array}
\label{Conn}
\end{equation}
\end{definition}

\subsection{Weil objects}

For a unital object $X$, we write $\Theta(X)$ for the set of
connections on $X.$ Then $\Theta$ is a functor from $\Unit(\mathcal{L})$ (therefore also from $\Asso(\mathcal{L}), \Com(  \mathcal{L})$) to $\Set$.

\begin{theorem}
$\Theta$ is a representable functor from  $\Unit(\mathcal{L}), \Asso(\mathcal{L}), \Com(\mathcal{L})$ to $\Set$. The representing objects, the \emph{Weil objects}, are defined as objects in $\mathcal{S}$ (and as $L_{0}$ modules) by
\begin{equation} \label{=repr}
M_{L}    \sim I\oplus DL_{0}^{\ast}, \,\,\,
A_{L}    =T_{1}(  M_{L})  \sim T(  DL_{0}^{\ast}), \,\,\,
C_{L}   =S_{1}(  M_{L})  \sim S(  DL_{0}^{\ast})
\end{equation}
\end{theorem}
For any connection $\theta$ on an object $X$ in any of the three categories, write $c_{\theta}$ for the corresponding ``Chern-Weil'' morphism from $M_{L}, A_{L}, C_{L}$ to $X$. We have the following analog of \cite[th\'eor\`eme 3]{car1}.

\begin{theorem}
For any two connections $\theta, \theta^{\prime}$ on $X$, the Chern-Weil morphisms $c_{\theta}$, $c_{\theta^{\prime}}$ are homotopic.
\end{theorem}

Following again \cite{car1}, we define the \emph{equivariant cohomology functor} as $H_{\rm{eq}}=H_{C_{L}}$. Then we have the following result which generalizes the computation of the
equivariant cohomology of a point.

\begin{proposition}
\label{point} $H_{\rm{eq}}(I)  =H(C_{L})  \approx
S(  L_{0}^{\ast})  ^{L_{0}}$
\end{proposition}

\subsection{Main results}

Our main results are three quasi-isomorphism theorems. The first one
requires only rigidity for $L_{0}$.

\begin{theorem}
\label{non-commutative} The map $A_{L}\rightarrow C_{L}$ is a quasi-isomorphism.
\end{theorem}

Assume now that $L_{0}$ is a \emph{quadratic} Lie algebra, i.e., there is an
$L_{0}$-module isomorphism $L_{0}^{\ast}\approx L_{0}$. Then the unital Weil
object $M_{L}$ is naturally a \emph{unital Lie\ algebra}, i.e. it admits a
distinguished Lie algebra morphism from the trivial Lie algebra $I$. In fact
one has the following central extension of Lie algebras:
$0\rightarrow I\rightarrow M_{L}\rightarrow L\rightarrow0$.
Now for any unital Lie algebra one can define its \emph{unital enveloping
algebra}, the universal object corresponding to the forgetful functor
from associative algebras to unital Lie algebras. We write $E_{L}$ for the unital enveloping algebra of $M_{L}$. We obtain a natural algebra morphism $A_{L}\rightarrow E_{L}$ :
\begin{theorem}
\label{quantum} The map $A_{L}\rightarrow E_{L}$ is a quasi-isomorphism.
\end{theorem}

Following \cite{am2} we can define noncommutative and quantized equivariant
cohomology (the latter for quadratic Lie algebras) : $\mathbf{H}_{\rm{eq}}=H_{A_{L}}\text{ and }\mathcal{H}_{\rm{eq}}=H_{E_{L}}$.

\begin{corollary}
For any algebra $B$ in $\Asso(  \mathcal{L})$, the algebras
$H_{\rm eq}(B)$, $\mathbf{H}_{\rm{eq}}(B)$ and $\mathcal{H}_{\rm{eq}}(B)$
are isomorphic.
\end{corollary}

\begin{proposition}\label{qChev}
$H(  E_{L})  \approx U(  L_{0})  ^{L_{0}}$.
\end{proposition}
Combining Propositions \ref{point}, \ref{qChev} and Theorems \ref{quasi}, \ref{non-commutative}, \ref{quantum}, we obtain for $L_{0}$ quadratic ({\it cf} \cite {dufl})~:

\begin{corollary}
[Duflo isomorphism] The algebras $U(  L_{0})  ^{L_{0}}$ and $S(  L_{0})^{L_{0}}$ are isomorphic.
\end{corollary}

\end{document}